\documentclass[10pt]{article}
\usepackage{amssymb}
\usepackage[english,francais]{babel}
\usepackage{amsmath}
\usepackage{bm}

\newtheorem{theorem}{Theorem}[section]

\newtheorem{e-proposition}[theorem]{Proposition}

\newtheorem{e-definition}[theorem]{Definition\rm}

\newtheorem{theoreme}{Th\'eor\`eme}[section]

\newtheorem{definition}[theoreme]{D\'efinition\rm}

\def\og{\leavevmode\raise.3ex\hbox{$\scriptscriptstyle\langle\!\langle$~}}
\def\fg{\leavevmode\raise.3ex\hbox{~$\!\scriptscriptstyle\,\rangle\!\rangle$}}

\begin{document}

\date{\empty}
\title{The motion of a rigid body and a viscous fluid in a bounded domain in
presence of collisions}
\author{{\large Nikolai V. Chemetov$^1$, \v S\'arka Ne\v casov\'a$^2$} \\
%EndAName
\\
{\small $^1$ University of Lisbon, Portugal }\\
{\small nvchemetov@fc.ul.pt}\\
{\small $^2$ Institute of Mathematics, }\\
{\small \v Zitn\'a 25, 115 67 Praha 1, Czech Republic }\\
{\small matus@math.cas.cz}}
\maketitle

\begin{abstract}
\selectlanguage{english}{We consider the motion of a rigid body, governed by the Navier-Stokes
equations in a bounded domain. Navier's condition is prescribed on the
boundary of the body. We give the global in a time solvability result of
weak solution. The result permits a possibility of collisions of the body
with the boundary of the domain.}

\vskip 0.5\baselineskip

\selectlanguage{francais} \vskip 0.5\baselineskip \noindent \textbf{Le
mouvement d'un corps rigide avec les conditions aux limites de Navier en
pr\'esence de collisions. } Nous consid\'erons le mouvement d'un corps
rigide, r\'egi par les \'equations de Navier-Stokes dans un domaine born\'e.
L'\'etat de Navier est prescrit \`a la limite du corps. Nous donnons le
global dans le r\'esultat de solvabilit\'e de temps de solution faible. Ce
r\'esultat permet une possibilit\'e de collisions du corps avec le bord du
domaine.
\end{abstract}

% place in the next line the header (rubrique) chosen for your article,
% if you know it (you can also have 2, format : Header1/Header2
\centerline{} %\begin{frontmatter}

\selectlanguage{english}

% The [label1] can be suppressed if there is only one address for all authors

\selectlanguage{english}

\medskip %\begin{center}
%{\small Received *****; accepted after revision +++++\\
%Presented by £££££}
%\end{center}

%\end{frontmatter}

% now the Version française abrégée, if it exists
\selectlanguage{francais}

\section*{Version fran\c{c}aise abr\'eg\'ee}

% Text of your Version française abrégée here.
% Note you do not need to repeat here equations that you use in the
% main text - for example 'voir (3)' is quite acceptable.

%In the article we investigate the motion of a rigid body inside of a viscous
%incompressible fluid.
Dans l'article, nous \'etudions le mouvement d'un corps rigide \`{a}
l'int\'erieur d'un visqueux fluide incompressible.
%Let $\Omega $ be a bounded domain of $\mathbb{R}^{N}$ for $N=2$ or $3$.
Soit $\Omega$ un domaine born\'e de $\mathbb{R}^{N}$, $N=2$ ou $3$.
%At the initial moment $t=0$ the body and the fluid occupy
%an open connected set ${S}_{0}\subset \Omega $ and
%the set $F_{0}=\Omega \backslash \overline{S_{0}},$ respectively.
A l'instant initial $t = 0$ le corps et le fluide occupent un ensemble
ouvert connexe ${S}_{0} \subset \Omega$ et l'ensemble $F_{0} = \Omega
\backslash \overline{S_{0}},$ respectivement. %At any time $t \in [0,T]$
%the motion of the fluid and the body is governed by
%system \eqref{is}-\eqref{eg33}.
A tout moment, $t \in [0,T]$ le mouvement du fluide et le corps est r\'egie
par ~ system \eqref{is}-\eqref{eg33}.

%The global existence of weak solution for this system with
%a non-slip boundary condition (Dirichlet's one) on boundaries of the body and the domain
%has been treated by many
%mathematicians.
L'existence d'une solution globale faible pour ce syst\`{e}me avec une
condition aux limites non-glissement (celui de Dirichlet) sur les limites du
corps et le domaine a \'et\'e trait\'ee par de nombreux math\'ematiciens.
%The investigations in the articles  \cite{GH}, \cite{HES},  \cite{HIL},
%\cite{STA2} has shown that the non-slip
%boundary condition gives a paradoxical result of no collisions between the
%body and the boundary of the domain.
Les enqu\^{e}tes dans les articles \cite{GH}, \cite{HES}, \cite{HIL}, \cite%
{STA2} a montr\'e que l'anti-d\'erapant condition limite donne un r\'esultat
paradoxal d'aucune collision entre le le corps et le bord du domaine.
%These results have demonstrate that a
%more accurate model is needed for the description of the
%motion of bodies in a viscous incompressible fluid.
Ces r\'esultats ont d\'emontr\'e qu'un mod\`{e}le plus pr\'ecis est
n\'ecessaire pour la description de la mouvement des corps dans un fluide
visqueux incompressible.

%In the articles \cite{GH2}, \cite{NP1},  \cite{PS111} have proposed to
%study system \eqref{is}-\eqref{eg33} by adding a slip boundary condition (Navier's one) instead of
%the non-slip boundary condition  on boundaries of the body and the domain.
Dans les articles \cite{GH2}, \cite{NP1}, \cite{PS111} ont propos\'e de syst%
\`{e}me \'etude \eqref{is}-\eqref{eg33} en ajoutant une condition de limite
de glissement (de Navier un) au lieu de la condition \`{a} la limite de
non-glissement sur les limites du corps et le domaine.
%In these articles it has been shown that the slip boundary condition
%cleans the no-collision paradox.
Dans ces articles il a \'et\'e montr\'e que la condition de limite de
glissement nettoie le paradoxe sans collision.

%In our article we close system \eqref{is}-\eqref{eg33} by adding
% Navier's \ boundary condition \eqref{eq3111} on the boundary of the body $S(t)$
%and  Dirichlet's boundary condition \eqref{2b} on the boundary of the domain $\Omega$.
Dans notre article, nous fermons syst\`{e}me \eqref{is}-\eqref{eg33} en
ajoutant ~ La condition aux limites de Navier \eqref{eq3111} sur la limite
du corps $S(t)$ et l'\'etat limite de Dirichlet \eqref{2b} sur la fronti\`{e}%
re du domaine $\Omega$.
%With these boundary conditions in Theorem   \ref{theorem} we have shown the existence
% of weak solution for problem \eqref{is}-\eqref{2b}.
Avec ces conditions aux limites dans le th\'eor\`{e}me \ref{theorem} nous
avons montr\'e l'existence ~ de solution faible pour le probl\`{e}me %
\eqref{is}-\eqref{2b}. %Due to the low regularity of the boundaries
%$\partial \Omega \in C^{0,1}$, $\partial {S}_{0}\in C^{2}$
%and that the external force
% $\mathbf{g}\in L^{2}((0,T); (LD_{0}^{2}(\Omega ))^{\ast })$,
%the contacts of the body and the
%boundary of the domain are available in Theorem \ref{theorem}.
En raison de la faible r\'egularit\'e des limites ~$\partial \Omega \in
C^{0,1}$, $\partial {S}_{0}\in C^{2}$ et que la force ext\'erieure $\mathbf{g%
}\in L^{2}((0,T); (LD_{0}^{2}(\Omega ))^{\ast })$, les contacts du corps et
la fronti\`{e}re du domaine sont disponibles dans le th\'eor\`{e}me \ref%
{theorem}.

%To show Theorem \ref{theorem}, we consider the approximate problem
%(\ref{a1})-(\ref{a3}) to system
%\eqref{is}-\eqref{2b}, using the idea that the "body+fluid"  can be approximated by a
%\textit{non-homogeneous} fluid, having different values of viscosity in
%three zones: approximation of "body", approximation of Navier's boundary
%condition \eqref{eq3111} and "fluid" \ zone.
Pour montrer le th\'eor\`{e}me \ref{theorem}, nous consid\'erons le probl%
\`{e}me approximative ~(\ref{a1})-(\ref{a3}) au syst\`{e}me \eqref{is}-%
\eqref{2b}, en utilisant l'id\'ee que le "corps+fluide" peut \^{e}tre
approch\'ee par une fluide non homog\`{e}ne, ayant des valeurs diff\'erentes
de viscosit\'e en trois zones: rapprochement des "corps", approximation de
la limite de Navier \'etat \eqref{eq3111} et la zone "fluide".
%Then in problem (\ref{a1})-(\ref{a3}) we pass on limits with respect of
%the parameters $\varepsilon ,\delta$.
Puis, en probl\`{e}me (\ref{a1})-(\ref{a3}) que nous passons sur les limites
\`{a} l'\'egard de les param\`{e}tres $\varepsilon ,\delta$.
%In the limit on $\varepsilon \rightarrow 0$ we obtain
%the motion of the rigid body in a viscous fluid.
Dans la limite $\varepsilon \rightarrow 0$ nous obtenons le mouvement du
corps rigide dans un fluide visqueux.
%Then in the second limit on $\delta \rightarrow 0$, we obtain the motion\
%of the  body already with Navier's boundary condition \eqref{eq3111}.
Puis, dans la deuxi\`{e}me limite $\delta \rightarrow 0$, nous obtenons le
mouvement du corps d\'ej\`{a} avec la condition limite de Navier %
\eqref{eq3111}.
%Since the demonstration of Theorem \ref{theorem} is a quite lengthy and technical
%one for details we refer to \cite{cnn}.
Depuis la d\'emonstration du th\'eor\`{e}me \ref{theorem} est une assez
longue et technique l'un pour plus de d\'etails, nous renvoyons \`{a} \cite%
{cnn}.

\selectlanguage{english} % main text

\section{Presentation of the problem}

\label{ms}

We investigate the motion of a rigid body inside of a viscous incompressible
fluid. Let $\Omega $ be a bounded domain of $\mathbb{R}^{N}$ for $N=2$ or $3$%
. At the initial moment $t=0$ the body and the fluid occupy an open
connected set ${S}_{0}\subset \Omega $ and the set $F_{0}=\Omega \backslash
\overline{S_{0}},$ respectively. The motion of any point $\mathbf{y}%
=(y_{1},..,y_{N})^{T}\in S_{0}$ is described by a preserving orientation
isometry
\begin{equation}
\mathbf{A}(t,\mathbf{y})=\mathbf{q}(t)+\mathbb{Q}(t)(\mathbf{y}-\mathbf{q}%
(0)), \qquad t\in \lbrack 0,T],  \label{is}
\end{equation}%
where $\mathbf{q}=\mathbf{q}(t)$ is the body mass center and $\mathbb{Q}=%
\mathbb{Q}(t)$ is the rotation matrix, such that $\mathbb{Q}(t)\mathbb{Q}%
(t)^{T}=\mathbb{I},$ $\mathbb{Q}(0)=\mathbb{I}$\ with $\mathbb{I}$ being the
identity matrix. Hence the body and the fluid occupy the sets $\
S(t)=A(t,S_{0})$ and $F(t)=\Omega \backslash \overline{S(t)}$ at any time $%
t. $ The velocity of the body is related with the isometry $\mathbf{A}$ by
\begin{equation}
\mathbf{u}=\mathbf{q}^{\prime }(t)+\mathbb{P}(t)(\mathbf{x}-\mathbf{q}%
(t))\qquad \text{for}\ \ \mathbf{x}\in S(t),  \label{comp}
\end{equation}%
where a matrix $\mathbb{P}(t)$ fulfills $\frac{d \mathbb{Q}} {dt} \mathbb{Q}%
^{T}=\mathbb{P}$, such that there exists a vector $\bm{\omega}=\bm{\omega}%
(t)\in \,\mathbb{R}^{N},$ satisfying $\mathbb{P}(t)\mathbf{x}=\bm{\omega}%
(t)\times \mathbf{x},$ $\ \forall \mathbf{x}\in \mathbb{R}^{N}.$

The motion of the fluid and the body is governed by the following system%
\begin{align}
\partial _{t}\rho +(\mathbf{u}\cdot \nabla )\rho & =0,\qquad \mathrm{div}%
\mathbf{u}=0,\qquad \quad \text{for}\;\mathbf{x}\in F(t),  \notag \\
\rho (\partial _{t}\mathbf{u}+(\mathbf{u}\cdot \nabla )\mathbf{u})& =\mathrm{%
div}P+\mathbf{g},  \notag \\
m\mathbf{q}^{\prime \prime }& =\int_{\partial S(t)}P\mathsf{n}\,d\mathbf{x}%
+\int_{S(t)}\mathbf{g}\,d\mathbf{x},\qquad \quad \text{for}\;\mathbf{x}\in
S(t),  \notag \\
\frac{d(\mathbb{J}\bm{\omega})}{dt}& =\int_{\partial S(t)}(\mathbf{x}-%
\mathbf{q}(t))\times P\mathsf{n}\,d\mathbf{x}+\int_{S(t)}(\mathbf{x}-\mathbf{%
q}(t))\times \mathbf{g}\,d\mathbf{x}  \label{system1}
\end{align}%
with the initial conditions%
\begin{equation}
S=S_{0},\quad \rho =\rho _{0},\quad \mathbf{u}=\mathbf{u}_{0}\quad \text{at }%
t=0.  \label{eg33}
\end{equation}%
Here $\ m=\int_{S(t)}\rho \,d\mathbf{x}-$ the mass of the body; $\rho $ is
the density in the body $S(t)$ and in the fluid $F(t)$; $P=-pI+2\mu _{f}\,%
\mathbb{D}\mathbf{u}$ and $\mathbb{D}\mathbf{u}=\frac{1}{2}\{\nabla \mathbf{u%
}+\left( \nabla \mathbf{u}\right) ^{T}\}$ -the stress and the
deformation-rate tensors of the fluid; $p$ -the fluid pressure; $\mu _{f}>0$
- the constant viscosity of the fluid; $\mathsf{n}$ -the unit outward normal
to $\partial S(t)$; $\mathbb{J}=\int_{S(t)}\rho (|\mathbf{x}-\mathbf{q}%
(t)|^{2}\mathbb{I}-(\mathbf{x}-\mathbf{q}(t))\otimes (\mathbf{x}-\mathbf{q}%
(t)))\,d\mathbf{x}$ -the matrix of the inertia moments of the body; $\mathbf{%
g}$-an external force.

The global existence of weak solution has been treated by many
mathematicians: Hoffmann, Starovoitov, Conca, San Mart\'in, Tucsnak,
Feireisl, Ne\v casov\'a, Hillairet, Bost, Cottet, Maitre; Desjardins,
Esteban, Gunzburger, Lee, Seregin, Takahashi and etc.. All of these authors
have considered non-slip boundary condition on boundaries of the body and
the domain, but this boundary condition \ gives a paradoxical result of no
collisions between the body and the boundary of the domain: Hesla \cite{HES}%
, Hillairet \cite{HIL}, Starovoitov \cite{STA2}. \ In the articles \cite{GH}%
, \cite{SST}, \cite{STA2} the autors have studied the question of possible
collisions with respect of the regularity of velocity and the regularity of
boundaries. \ For instance, in \cite{GH} G\'erard-Varet, Hillairet have
demonstrated that under $C^{1,\alpha }$-boundaries the collision is possible
in finite time if and only if $\alpha <1/2.$ These mentioned results have
demonstrated that a more accurate model is needed for the description of the
motion of bodies in a viscous incompressible fluid.

Neustupa, Penel \cite{NP1} have investigated a prescribed collision of a
ball with a wall, when the slippage is allowed on the boundaries of the ball
and of the wall. The slippage is prescribed by Navier's boundary condition,
having only the continuity of velocity field just in the normal component.
This pioneer result \cite{NP1} have shown that the slip boundary condition
cleans the no-collision paradox. Recently G\'erard-Varet, Hillairet \cite%
{GH2} have proved a local-in-time existence result: up to collisions. The
motion of a single body, moved in the whole space $\mathbb{R}^{3},$ have
considered in \cite{PS111}. The free fall of a ball above a wall with the
slippage, prescribed on the boundaries, have been studied in \cite{GHC},
where it was shown that the ball touches the boundary of the wall in a
finite time.

In this article we close system \eqref{system1} by adding Navier's \
boundary condition
\begin{equation}
\mathbf{u}_{s}\cdot \mathsf{n}=\mathbf{u}_{f}\cdot \mathsf{n},\qquad (P_{f}%
\mathsf{n}+\gamma (\mathbf{u}_{f}-\mathbf{u}_{s}))\cdot \mathsf{s}=0\quad
\text{on }\partial S(t),  \label{eq3111}
\end{equation}%
and Dirichlet's boundary condition%
\begin{equation}
\mathbf{u}=0\quad \text{on }\partial \Omega .  \label{2b}
\end{equation}%
Here $\mathbf{u}_{s}$ and $\mathbf{u}_{f}$ are the trace values of the
velocity $\mathbf{u}$ on $\partial S(t)$ from the rigid\ side $S(t)$ and
from the fluid \ side $F(t),$ respectively; $\mathsf{n}$ and $\mathsf{s}$
are the external normal and arbitrary tangent vector to $\partial S(t);$ the
constant $\gamma >0$ is the friction coefficient of $\partial S_{0}.$

\section{Weak solution of system \eqref{is}-\eqref{2b} and the main result}

\label{ms2}

To introduce the concept of weak solution for system \eqref{is}-\eqref{2b},
let us define some spaces of functions%
\begin{eqnarray*}
V^{0,2}(\Omega ) &=&\{\mathbf{v}\in L^{2}(\Omega ):\,\mbox{div }\mathbf{v}%
=0\quad \text{in}\;\mathcal{D}^{\prime }(\Omega ),\quad \mathbf{v}\cdot
\mathsf{n}=0\;\ \text{ in}\;H^{-1/2}(\partial \Omega )\}, \\
BD_{0}(\Omega ) &=&\left\{ \mathbf{v}\in L^{1}(\Omega ):\,\mathbb{D}\mathbf{v%
}\in \mathcal{M}(\Omega ),\quad \mathbf{v}=0\quad \text{on}\;\partial \Omega
\right\} ,
\end{eqnarray*}%
where $\mathsf{n}$ is the unit normal to the boundary $\partial \Omega $ of
the domain $\Omega $ and $\mathcal{M}(\Omega )$ is the space of bounded
Radon measures. \ Let $S$\ be an open connected subset of $\Omega .$ \ We
consider the space%
\begin{equation*}
KB(S)=%
\begin{array}{c}
\left\{ \mathbf{v}\in BD_{0}(\Omega ):\,\mathbb{D}\mathbf{v}\in L^{2}(\Omega
\backslash \overline{S}),\quad \mathbb{D}\mathbf{v}=0\quad \text{a.e. on }%
S,\right. \\
\left. \mbox{div}\mathbf{v}=0\quad \mbox{
in }\ \mathcal{D}^{\prime }(\Omega )\right\}%
\end{array}%
.
\end{equation*}

\begin{definition}
\label{definition} The triple $\left\{ \mathbf{A},\rho ,\mathbf{u}\right\} $
is a weak solution of system \eqref{is}-\eqref{2b}, if the following three
conditions are fulfilled:

1) The function $\mathbf{A}(t,\cdot ):\mathbb{R}^{N}\rightarrow \mathbb{R}%
^{N}$ \ is a preserving orientation isometry \eqref{is}, such that\ the
functions $\mathbf{q},$ $\mathbb{Q}$ are absolutely continuous on $[0,T]$. \
The isometry $\mathbf{A}$ is compatible with the rigid body velocity %
\eqref{comp} on $S(t)$ and defines a time dependent set $S(t)=\mathbf{A}%
(t,S_{0})$;

2) The function $\rho \in L^{\infty }((0,T)\times \Omega )$ satisfies the
integral equality%
\begin{equation}
\int_{0}^{T}\int_{\Omega }\rho (\xi _{t}+(\mathbf{u}\cdot \nabla )\xi )\,dtd%
\mathbf{x}=-\int_{\Omega }\rho _{0}\xi (0,\cdot )\,d\mathbf{x}  \label{eq5}
\end{equation}%
for any $\,\,\xi \in C^{1}([0,T]\times \overline{\Omega }),\quad \xi
(T,\cdot )=0;$

3) The function $\mathbf{u}\in L^{2}(0,T;KB(S(t)))\cap L^{\infty
}(0,T;V^{0,2}(\Omega ))$ satisfies the integral equality%
\begin{align}
& \int_{0}^{T}\int_{\Omega \backslash \partial S(t)}\{\rho \mathbf{u}%
\boldsymbol{\psi }_{t}+\rho (\mathbf{u}\otimes \mathbf{u}):\mathbb{D}%
\boldsymbol{\psi }-2\mu _{f}\,\mathbb{D}\mathbf{u}:\mathbb{D}\boldsymbol{%
\psi }\,+\mathbf{g}\boldsymbol{\psi }\}d\mathbf{x}dt  \notag \\
& =-\int_{\Omega }\rho _{0}\mathbf{u}_{0}\boldsymbol{\psi }(0,\cdot )\,d%
\mathbf{x}+\int_{0}^{T}\left\{ \int_{\partial S(t)}\beta (\mathbf{u}_{s}-%
\mathbf{u}_{f})(\boldsymbol{\psi }_{s}-\boldsymbol{\psi }_{f})\,d\mathbf{x}%
\right\} dt  \label{eq4}
\end{align}%
for any $\boldsymbol{\psi }\in L^{2}(0,T;KB(S(t))),$ such that $\boldsymbol{%
\psi }_{t}\in L^{2}(0,T;L^{2}(\Omega \backslash \partial S(t)))~$ and $%
\boldsymbol{\psi }(T,\cdot )=0.$ Here\ we denote the trace values of $%
\mathbf{u},$ $\boldsymbol{\psi }$ \ on $\partial S(t)$\ \ from \ the \textit{%
rigid} side $S(t)$ and the fluid side $F(t)$ by $\mathbf{u}_{s}(t,\mathbf{%
\cdot }),$ $\boldsymbol{\psi }_{s}(t,\mathbf{\cdot })$ and $\mathbf{u}_{f}(t,%
\mathbf{\cdot }),$ $\boldsymbol{\psi }_{f}(t,\mathbf{\cdot }),$ respectively.
\end{definition}

Our main result is the following theorem.

\begin{theorem}
\label{theorem} We assume that $S_{0}\subset \Omega ,$ such that $%
dist[S_{0},\Omega ]>0.$ We admit that the boundaries $\partial \Omega \in
C^{0,1}$ and $\partial {S}_{0}\in C^{2}.$ Let
\begin{equation*}
\rho _{0}(\mathbf{x})=\left\{
\begin{array}{ll}
\rho _{s}(\mathbf{x})\geqslant const>0, & \quad \mathbf{x}\in S_{0}; \\
\rho _{f}=const>0, & \quad \mathbf{x}\in F_{0},%
\end{array}%
\right. \qquad \rho _{s}\in L^{\infty }(S_{0}),
\end{equation*}%
\begin{equation}
\qquad \mathbf{u}_{0}\in V^{0,2}(\Omega ),\quad \mathbb{D}\mathbf{u}%
_{0}=0\quad \text{in \ }\mathcal{D}^{\prime }(S_{0}),\quad \mathbf{g}\in
L^{2}((0,T);(LD_{0}^{2}(\Omega ))^{\ast }).  \label{data}
\end{equation}%
Then system \eqref{is}-\eqref{2b}\ \ possesses a weak solution $\left\{
\mathbf{A},\rho ,\mathbf{u}\right\} ,$ such that the isometry $\mathbf{A}%
(t,\cdot )$ is Lipschitz continuous with respect to $t\in \lbrack 0,T],$%
\begin{equation}
\rho (t,\mathbf{x})=\left\{
\begin{array}{ll}
\rho _{s}(\mathbf{A}^{-1}(t,\mathbf{x})), & \mathbf{x}\in S(t); \\
\rho _{f}, & \mathbf{x}\in F(t),%
\end{array}%
\right. \quad \text{for a.e. }t\in (0,T),  \label{umm}
\end{equation}%
$\mathbf{u}\in C_{\mathrm{weak}}(0,T;V^{0,2}(\Omega ))$ and the following
energy inequality holds
\begin{align}
\frac{1}{2}\int_{\Omega }\rho |\mathbf{u}|^{2}(r)\ d\mathbf{x}&
+\int_{0}^{r}\left\{ \int_{F(t)}2\mu _{f}\,|\mathbb{D}\,\mathbf{u}|^{2}\,\,d%
\mathbf{x}+\int_{\partial S(t)}\beta |\mathbf{u}_{f}-\mathbf{u}_{s}|^{2}\ d%
\mathbf{x}\right\} dt  \notag \\
& \leqslant \frac{1}{2}\int_{\Omega }\rho _{0}|\mathbf{u}_{0}|^{2}\ d\mathbf{%
x+}\int_{0}^{r}<\mathbf{g,u}>\ dt\quad \text{for a.e. }\,r\in (0,T).
\label{energy}
\end{align}
\end{theorem}

Let us point that in \cite{GHC} it has been shown that the ball never
touches the boundary of the wall for mixed boundary conditions \eqref{eq3111}%
, \eqref{2b}. Nevertheless of the result \cite{GHC}, the contacts of the
body and the boundary of the domain are available in Theorem \ref{theorem},
due to the low regularity of the boundaries $\partial \Omega \in C^{0,1}$, $%
\partial {S}_{0}\in C^{2}$. And moreover $\mathbf{g}\in L^{2}((0,T);
(LD_{0}^{2}(\Omega ))^{\ast })$, we refer to the example constructed by
Starovoitov \cite{STA2}. In order to create a collision of the body with the
boundary of the domain (in the case of non-slip conditions on the boundaries
$\partial \Omega $ and $\partial {S}_{0}$), Starovoitov adds an external
force from $H^{-1}$-space.

\section{\textbf{Sketch of the proof of Theorem} \ }

First we introduce an approximate scheme to system \eqref{is}-\eqref{2b},
using the idea that the "body+fluid" can be approximated by a \textit{%
non-homogeneous} fluid, having different values of viscosity in three zones:
approximation of "body", approximation of Navier's boundary condition %
\eqref{eq3111} and "fluid" \ zone.

To construct such approximation problem we fix the following notations. For
an open connected set $S\subset \mathbb{R}^{N}$, we define $\mathrm{dist}[%
\mathbf{x},S]=\inf_{\mathbf{y}\in S}|\mathbf{x}-\mathbf{y}|,$ $\ d_{S}(%
\mathbf{x})=\mathrm{dist}[\mathbf{x},{\,\mathbb{R}^{N}\setminus S}]-\mathrm{%
dist}[\mathbf{x},S]$\ for any\ $\mathbf{x}\in \mathbb{R}^{N}$ and $\
[S]_{\delta }=d_{S}^{-1}((\delta ,+\infty ))$-the $\delta -$kernel of $S$, $%
\ ]S[_{\delta }=d_{S}^{-1}((-\delta ,+\infty ))$-the $\delta -$neighborhood
of $S$. The characteristic functions of the sets \ $S_{0}$, $[S_{0}]_{\delta
}$ and $]S_{0}[_{\delta }\backslash \overline{S_{0}}$ are denoted by $%
\varphi _{0}$, $\zeta _{0}^{\delta }$ \ and $\chi _{0}^{\delta }, $ which
are\ defined in $\mathbb{R}^{N}.$

The approximation problem to system \eqref{is}-\eqref{2b} consists from the
linear transport equations%
\begin{eqnarray}
\partial _{t}\rho +(\mathbf{u}\cdot \nabla )\rho  &=&0,\qquad \partial
_{t}\zeta +(\overline{\mathbf{u}}^{\delta }\cdot \nabla )\zeta =0\qquad
\text{in}\;(0,T)\times \mathbb{R}^{N},  \notag \\
\rho (0) &=&\rho _{0}^{\varepsilon },\qquad \zeta (0)=\zeta _{0}^{\delta
}\qquad \text{in}\;\mathbb{R}^{N},  \label{a1}
\end{eqnarray}%
and the momentum equation%
\begin{eqnarray}
\rho (\partial _{t}\mathbf{u}+(\mathbf{u}\cdot \nabla )\mathbf{u}) &=&%
\mathrm{div}P+\mathbf{g},\quad \mathrm{div}\mathbf{u}=0\quad \text{in}%
\;(0,T)\times \Omega ,  \notag \\
\mathbf{u}(0) &=&\mathbf{u}_{0}\quad \text{in}\;\Omega ,  \label{a2}
\end{eqnarray}%
where $\rho _{0}^{\varepsilon }=\varepsilon \chi _{0}^{\delta }+\rho
_{f}\theta _{0}^{\delta }+\rho _{s}\varphi _{0},$ \ $\theta _{0}^{\delta
}=1-(\varphi _{0}+\chi _{0})$ \ \ and%
\begin{equation}
P=-p+\mu _{\varepsilon ,\delta }\mathbb{D}\mathbf{u},\qquad \mu
_{\varepsilon ,\delta }=\frac{1}{\varepsilon }\varphi +2\delta \beta \chi
+2\mu _{f}\theta ,\qquad \theta =1-(\varphi +\chi )  \label{a3}
\end{equation}%
and $\varphi $, $\chi $\ \ are the characteristic functions, defined in $%
\mathbb{R}^{N}$, of the sets $S(\varphi )=]S(\zeta )[_{\delta }$, \ $S(\chi
)=]S(\zeta )[_{2\delta }\backslash \overline{]S(\zeta )[_{\delta }}$. The
function $\overline{\mathbf{u}}^{\delta }$ is the standard mollification of $%
\mathbf{u}$ on the parameter $\delta .$

In relation (\ref{a3}) the $\varepsilon -$dependence of the
\textquotedblleft viscosity\textquotedblright\ $\mu _{\varepsilon ,\delta }$
is a penalization, introduced in \cite{HOST}, where the rigid body is
replaced by a fluid, having high viscosity value. The $\delta -$dependence
of $\mu _{\varepsilon ,\delta }$ defines a mixture region between the fluid
and the "body", which approximates the jump boundary term on $\partial S(t)$
in (\ref{eq4}). The solvability of this approximation problem (\ref{a1})-(%
\ref{a3}) can be shown by a fixed point argument, Galerkin's method and
theoretical results for transport equations (see \cite{ant}, \cite{boyer}),
\cite{Temam}).

Next in the approximation problem we have to pass on limits with respect of
the parameters $\varepsilon ,\delta$. These limits are based on the results
for the transport equations \cite{boyer}.

\begin{itemize}
\item The first limit on $\varepsilon \rightarrow 0$ is related with a
so-called \textit{"solidification"} \ procedure in the zone of the
non-homogeneous fluid, corresponding to the "body". This limit can be
treated as in \cite{HOST}, \cite{SST}. In the limit we obtain the motion of
the rigid body in a viscous fluid;

\item In the second limit on $\delta \rightarrow 0$, we obtain the motion\
of the body already with a prescribed Navier's boundary condition. Firstly
we need to construct an appropriate set of test functions, depending on $%
\delta .$ Then, using embedding results in cusp domains, we show that the
viscous term $2\delta \beta \chi $ converges to the jump boundary term on $%
\partial S(t)$ in (\ref{eq4}). The embedding results allow also to apply a
compactness result in the convective term of (\ref{a2}) by using the
approach of Proposition 6.1 in \cite{SST}.\ \bigskip
\end{itemize}

The demonstration of Theorem \ref{theorem} is a quite lengthy and technical
one. The details can be found in \cite{cnn}.

% etc, etc

% \appendix

% \section{}
% \label{}

% The Acknowledgements are an un-numbered section
%\section*{Acknowledgements}
% Acknowledgements text here

\end{document}